# On the Identifiability of Latent Models for Dependent Data

Stéphane Guerrier† and Roberto Molinari‡

**Abstract:** The condition of parameter identifiability is essential for the consistency of all estimators and is often challenging to prove. As a consequence, this condition is often assumed for simplicity although this may not be straightforward to assume for a variety of model settings. In this paper we deal with a particular class of models that we refer to as "latent" models which can be defined as models made by the sum of underlying models, such as a variety of linear state-space models for time series. These models are of great importance in many fields, from ecology to engineering, and in this paper we prove the identifiability of a wide class of (second-order stationary) latent time series and spatial models and discuss what this implies for some extremum estimators, thereby reducing the conditions for their consistency to some very basic regularity conditions. Finally, a specific focus is given to the Generalized Method of Wavelet Moments estimator which is also able to estimate intrinsically second-order stationary models.

## 1. Introduction

The estimation of parametric models heavily relies on the condition that the model parameters are *identifiable* based on a specific estimation method. Supposing we have a parametric model denoted as $F_{\boldsymbol{\theta}}$, where $\boldsymbol{\theta}$ represents the $p$-dimensional parameter vector of interest, this condition implies that there is a unique maximizer (or minimizer) to the objective function which defines the estimator and that this unique solution is the true parameter value $\boldsymbol{\theta}_0$. More formally, let us consider an extremum estimator (see e.g. Newey and McFadden, 1994) defined as follows

$$\hat{\boldsymbol{\theta}} = \operatorname*{argmax}_{\boldsymbol{\theta} \in \boldsymbol{\Theta}} Q_n(\boldsymbol{\theta}),$$

where $\boldsymbol{\Theta}$ is the compact set of possible parameter values and $Q_n(\boldsymbol{\theta})$ is an objective function based on a sample size $n$ which we assume converges uniformly in probability to $Q(\boldsymbol{\theta})$. Based on this, we can define identifiability as being the case where, for all parameter vectors $\boldsymbol{\theta} \neq \boldsymbol{\theta}_0$, we have that

$$\sup_{\boldsymbol{\theta} \in \boldsymbol{\Theta} \setminus \mathcal{N}} Q(\boldsymbol{\theta}) < Q(\boldsymbol{\theta}_0),$$

where $\mathcal{N}$ is any open subset of $\boldsymbol{\Theta}$ containing $\boldsymbol{\theta}_0$. Consistency and asymptotic normality of all estimators depend on this condition which is often assumed for simplicity or proven on a case-by-case basis. The importance of this property was for example emphasized in Newey and McFadden (1994) where they state that, although difficult to prove in general, "*it is important to check that it is not a vacuous assumption whenever possible, by showing identification in some special cases*".

In this paper we prove and discuss the identifiability of a specific class of second-order stationary and regularly spaced time series and spatial models and refer to this class as "*latent*" models that consist in the sum of different underlying models which cannot directly be observed. These models are greatly used if we consider a broad class of linear state-space and structural time series models (see e.g. Durbin and Koopman, 2012) where the observed data can be explained by a series of components that are allowed to vary over time and have a direct interpretation. For example, in a wide variety of scientific domains many phenomena are explained by unobserved components and are also measured using specific devices which are unable to perfectly recover the values of the phenomenon of interest. A simple example of this is the case where a physical phenomenon, such as the behaviour of a chemical or biological substance, is measured in time but the measuring device is characterized by a stochastic error process which can also vary over time and/or location. The resulting measurement is therefore a combination of the realizations of a certain phenomenon plus a measurement error which the scientist would need to recover in order to have correct interpretations (see e.g. Buonaccorsi, 2010). Other examples where the interest lies in the unobserved components can be found in economics, ecology or biology where the presence of coexisting stochastic processes for different factors is quite common (see e.g. Harvey and Koopman, 1993; Cedersund and Roll, 2009; Zhang et al., 2010).

Despite being frequently used in many contexts, the identifiability of these models is generally assumed, even though this assumption may appear very strong for a variety of latent models. A class of models for which identifiability results do exist are the AutoRegressive Moving-Average (ARMA) models that, as emphasized by Granger and Morris (1976), can be obtained from a sum of latent autoregressive and white noise processes. However, it is unclear if it is always possible to recover the parameters of the latent models through the parameters of the corresponding ARMA as emphasized for example in Harvey and Koopman (1993). Moreover, if this were to be the case, these results could not be used when considering latent models composed by underlying models that are different from those that usually deliver an ARMA model. Therefore, although being widely used in practice, the question of whether these models are generally identifiable remains largely open thereby casting doubts on the consistency of the estimators which are used to recover their parameters.

Considering the above motivations, this paper proves the identifiability of a class of latent models for dependent data settings through their covariance function which is essential for any estimator to obtain the parameters of these models. Different latent models for time series as well as spatial data are investigated, consequently allowing to use this flexible class of models where a sum of simple spatial models can deliver a class of models in some way akin to ARMA models in the time series setting. Based on these findings, this pa-per also discusses what these imply for some extremum estimators, namely the Maximum Likelihood Estimators (MLE), Generalized Method of Moments estimators (GMM) as well as the recently proposed Generalized Method of Wavelet Moments (GMWM) estimators (see Guerrier et al., 2013). A specific discussion is developed for the latter class of estimators which can estimate models that include *intrinsically* stationary processes (i.e. non-stationary processes with stationary backward differences) such as drifts and random walks. The structure of the paper is therefore as follows. In Section 2 we present several results related to the identifiability of latent time series and spatial models while in Section 3 we discuss what these results imply for the mentioned extremum estimators, with a particular focus on GMWM estimators and a set of open-problems for future research, as well as a numerical example to support these findings. Finally, Section 4 concludes.

## 2. Injection of Covariance Models

The estimation of parameters of second-order stationary models for dependent data relies (almost) entirely on the covariance structure of these models. For this reason, in this section we give some results related to the identifiability for different sets of latent time series and spatial models through their respective covariance functions. To do so, let us first define the basic time series models considered in this work as $\left(X_t^{(j)}\right)$ with $t = 1 \ldots, T$ ($T \in \mathbb{N}_+$) and $j$ indicating a specific model. With this notation, we have:

(T1) *White Noise* (WN) with parameter $\sigma^2 \in \mathbb{R}^+$. We denote this process as $\left(X_t^{(1)}\right)$.

(T2) *Quantization Noise* (QN) (or rounding error, see e.g. Papoulis, 1991) with parameter $Q^2 \in \mathbb{R}^+$. We denote this process as $\left(X_t^{(2)}\right)$.

(T3) *Drift* with parameter $\omega \in \mathbb{R}^+$. We denote this process as $\left(X_t^{(3)}\right)$.

(T4) *Random walk* (RW) with parameter $\gamma^2 \in \mathbb{R}^+$. We denote this process as $\left(X_t^{(4)}\right)$.

(T5) *Moving Average* (MA(1)) process with parameter $\varrho \in (-1, +1)$ and $\varsigma^2 \in \mathbb{R}^+$. We denote this process as $\left(X_t^{(5)}\right)$.

(T6) *Auto-Regressive* (AR(1)) process with parameters $\rho \in (-1, +1)$ and $\upsilon^2 \in \mathbb{R}^+$. We denote this process as $\left(X_t^{(j)}\right)$, $j = 6, \ldots, K$ with $K \in \mathbb{N}^+$, $6 \leq K < \infty$.

Processes (T1), (T2), (T5) and (T6) are stationary models based on their parameter definitions. Process (T2) is particularly useful in the field of engineering where measurements are often rounded therefore introducing an error process in the procedure. In general terms, this process can be represented as a linear combination of differences of standard uniform variables where $Q^2$ plays the role of scaling factor. Processes (T3) and (T4) are typically non-stationary processes, where process (T3) is a non-random linear function with slope $\omega$, and will be discussed in Section 3.1 since their covariance functions are not defined. Having defined the time series models of interest for this paper, let us now define the spatial models. With $i \in \mathcal{I} \subset \mathbb{N}_+^2$ indicating the spatial coordinates of the observations, we have:

(S1) *Exponential model* with parameters $\phi \in \mathbb{R}^+$ and $\sigma^2 \in \mathbb{R}^+$. We denote this process as $\left(Y_i^{(k)}\right)$, $k = 1, \ldots, N_1$ with $N_1 \in \mathbb{N}^+$, $1 \leq N_1 < \infty$.

(S2) *Gaussian model* with parameters $\phi \in \mathbb{R}^+$ and $\sigma^2 \in \mathbb{R}^+$. We denote this process as $\left(Y_i^{(k)}\right)$, $k = 1, \ldots, N_2$ with $N_2 \in \mathbb{N}^+$, $1 \leq N_2 < \infty$.

These two spatial models are frequently used to approximate and describe many spatial processes. Moreover, a latent model made by the sum of either $N_1$ (S1) processes or $N_2$ (S2) processes represents a flexible solution to describe a more complex spatial process as mentioned in the introduction.

Based on the above definitions, let us now define the latent models for time series as $W_t = \sum_{j=1}^{K} X_t^{(j)}$ and for spatial processes as $V_i = \sum_{j=1}^{N} Y_t^{(j)}$, where $(V_i)$ can either be the sum of $N_1$ (S1) processes or of $N_2$ (S2) processes. Considering these models, let us denote their covariance function as $\varphi_h(\boldsymbol{\theta})$, where $h$ is the lag/distance between observations and hence $\boldsymbol{\varphi}(\boldsymbol{\theta}) = [\varphi_h(\boldsymbol{\theta})]_{h=1,\ldots,H}$. By finally denoting $\mathbf{A}(\boldsymbol{\theta}) = \partial/\partial\boldsymbol{\theta}\, \boldsymbol{\varphi}(\boldsymbol{\theta})$, let us set the following conditions:

**(C1)** All processes composing a latent model are mutually independent.
**(C2)** $\mathbf{A}(\boldsymbol{\theta})$ exists and is continuous $\forall \boldsymbol{\theta} \in \boldsymbol{\Theta}$.

Condition (C1) allows us to have a more tractable and general problem for the identifiability of the parameters of a latent model. Indeed, introducing dependence between the processes entails an approach to proving identifiability which is necessarily on a case-by-case basis since there can be different combinations of dependence between processes and their parametric dependence needs to be specified as well. Nevertheless, the results in this section give a sound basis to obtain further identifiability results also in these cases. As for Condition (C2), this is needed to develop an expansion of the covariance function $\boldsymbol{\varphi}(\boldsymbol{\theta})$ in order to prove that this function is injective. In fact, an approach that simultaneously verifies the usual identifiability conditions is to understand if the Jacobian matrix $\partial/\partial\boldsymbol{\theta}\, \boldsymbol{\varphi}(\boldsymbol{\theta})$ is of full column rank as a consequence of the following MacLaurin expansion

$$\boldsymbol{\varphi}(\boldsymbol{\theta}_1) = \boldsymbol{\varphi}(\boldsymbol{\theta}_0) + \mathbf{A}(\boldsymbol{\theta}^*) \underbrace{(\boldsymbol{\theta}_1 - \boldsymbol{\theta}_0)}_{\mathbf{b}} \tag{2.1}$$

where $\boldsymbol{\theta}_0$ and $\boldsymbol{\theta}_1$ are two parameter vectors and $\|\boldsymbol{\theta}^* - \boldsymbol{\theta}_0\| \leq \|\boldsymbol{\theta}_1 - \boldsymbol{\theta}_0\|$. Indeed, if $\boldsymbol{\theta}_0 = \boldsymbol{\theta}_1$, then we automatically have that $\boldsymbol{\varphi}(\boldsymbol{\theta}_0) = \boldsymbol{\varphi}(\boldsymbol{\theta}_1)$ but if we have $\boldsymbol{\theta}_0 \neq \boldsymbol{\theta}_1$ then, if the matrix $\mathbf{A}(\boldsymbol{\theta})$ is full column rank, it means that only the vector $\mathbf{b} = \mathbf{0}$ can make $\boldsymbol{\varphi}(\boldsymbol{\theta}_0) = \boldsymbol{\varphi}(\boldsymbol{\theta}_1)$ implying that the only situation where this is possible is when $\boldsymbol{\theta}_0 = \boldsymbol{\theta}_1$.

Knowing this, let us now define the classes of latent models studied in this section:

**Model 1** $W_t = X_t^{(1)} + X_t^{(2)} + \sum_{i=6}^{G} X_t^{(i)}$.
**Model 2** $W_t = X_t^{(5)} + \sum_{i=6}^{G} X_t^{(i)}$.
**Model 3** $V_i = \sum_{j=1}^{N} Y_i^{(j)}$, with $Y_i^{(j)}$ representing either model (S1) or model (S2).

The above processes are all second-order stationary ones, with **Model 1** and **Model 2** representing latent time series and **Model 3** representing latent spatial models. Since the covariance function is well defined for second-order stationary processes, we start by considering the class of models included within **Model 1** which represents the sum of a (T1) process, a (T2) process and $K$ (T6) processes (with $K = G - 5 < \infty$). Models of this class are commonly used in various disciplines going from economics to engineering for sensor calibration. Therefore, denoting $\rho_i$ as the autoregressive parameter of the $i^{th}$ (T6) process, the following theorem considers the properties of its autocovariance function with respect to the parameters of **Model 1**.

**Theorem 2.1.** *Under Conditions (C1) and (C2), and assuming $\rho_i \neq 0$, $\forall i$ and $\rho_i < \rho_j$, with $1 \leq i < j \leq K$, we have that the covariance function of* **Model 1** *is injective.*

The proof of this theorem can be found in Appendix A.1. Theorem 2.1 therefore states that there is a unique mapping of the parameters $\boldsymbol{\theta}$ to the covariance function $\boldsymbol{\varphi}(\boldsymbol{\theta})$ for a very general class of second-order stationary latent models. It must be noted that process (T5) is not included in this class due to the fact that the structure of matrix $\mathbf{A}(\boldsymbol{\theta})$ is not clearly full rank if this process is included in a latent model with processes (T1) and/or (T2). Based on this, we can give the following corollary to Theorem 2.1.

**Corollary 2.1.** *Under the conditions of Theorem 2.1, we have that the covariance function of* **Model 2** *is injective.*

The proof of this corollary is a direct consequence of the proof of Theorem 2.1 where the columns of matrix $\mathbf{A}(\boldsymbol{\theta})$ containing the derivatives with respect to the parameters of processes (T1) and (T2) are replaced by those of the parameters of process (T5).

**Remark 2.1.** *From the results of Granger and Morris (1976) we know that the sum of (T1) and (T6) models delivers different kinds of ARMA models. However, as emphasized in the introduction, it is not clear that all combinations of (T1) and (T6) models enjoy a unique mapping to a specific ARMA model as highlighted by Harvey and Koopman (1993). A simple example where this* mapping is indeed satisfied is given by the sum of two (T6) processes which has a unique mapping to an ARMA(2,1) model (see Appendix A.2 for details).

The results of Theorem 2.1 and Corollary 2.1 state that the covariance function of a wide class of second-order stationary latent time series models is injective. Let us now consider the latent processes $(V_i)$ based on the sum of spatial models (S1) and (S2) and give the following theorem.

**Theorem 2.2.** *Under Condition (C1) and (C2), and assuming $\phi_i < \phi_k$, $\forall\, i < k$, , with $1 \leq i < k \leq N$, the covariance function of Model 3 is injective.*

The result of this theorem delivers a class of latent spatial models which in some way presents some similarities with ARMA models for time series. Indeed, as mentioned earlier, a sum of (T6) models delivers an ARMA model based on the results of Granger and Morris (1976) and it is possible to suppose that a sum of Exponential or Gaussian covariance models can deliver something similar in the spatial case. Hence, Theorem 2.2 states that this new class of spatial models is identifiable based on their covariance function.

## 3. Consistency of Extremum Estimators for Latent Models

The results in Section 2 are helpful to considerably reduce the conditions for consistency of various estimators with respect to a wide class of second-order stationary models. In this section we consider some estimators that are part of the class of extremum estimators defined in the introduction and state the conditions for their consistency. The asymptotic distribution of these estimators is not discussed since it would require additional model-specific conditions but, in any case, it often relies strongly on the results of consistency given below. Having emphasized this, let us now list the other conditions (aside from identifiability) necessary to achieve consistency of extremum estimators.

**(C4)** $\Theta$ is compact.
**(C5)** $Q(\boldsymbol{\theta})$ is continuous.
**(C6)** $\sup_{\boldsymbol{\theta} \in \Theta} |Q_n(\boldsymbol{\theta}) - Q(\boldsymbol{\theta})| \xrightarrow{\mathcal{P}} 0$.

These conditions are based on Theorem 2.1 in Newey and McFadden (1994) and are standard for the consistency of extremum estimators. Indeed, Conditions (C4) and (C5) are statements which fall under the usual regularity assumptions while Condition (C6) requires the sample objective function $Q_n(\boldsymbol{\theta})$ to converge uniformly in probability to the true function $Q(\boldsymbol{\theta})$ which is not necessarily strong if adding minor assumptions to the stationary processes considered in this paper. With this in mind, the following paragraphs discuss these conditions for different estimators in the light of the results of Section 2.

### 3.1. MLE and GMM Estimators

The MLE and GMM estimators are very popular members of the class of extremum estimators and, considering these two particular estimators, let us de-note $f(\boldsymbol{z}|\boldsymbol{\theta})$ as the density function of the latent time series or spatial model while $g(\boldsymbol{z}|\boldsymbol{\theta})$ and $\boldsymbol{W}$ denote the moment conditions and the weighting matrix for the GMM estimators respectively. Considering this, let us give the following conditions:

**(C7)** $F_{\boldsymbol{\theta}}$ is uniquely defined by its covariance structure $\boldsymbol{\varphi}(\boldsymbol{\theta})$ and
$\mathbb{E}[|\ln(f(\boldsymbol{z}|\boldsymbol{\theta}))|] < \infty, \forall \boldsymbol{\theta}$.

**(C8)** The moment condition $g(\boldsymbol{z}|\boldsymbol{\theta})$ depends on $\boldsymbol{\theta}$ uniquely through the covariance function $\boldsymbol{\varphi}(\boldsymbol{\theta})$ and the weighting matrix $\boldsymbol{W}$ is positive semi-definite.

These conditions are based on Lemmas 2.2 and 2.3 in Newey and McFadden (1994) and list the additional properties needed for $Q(\boldsymbol{\theta})$ to be uniquely maximised at the true parameter value (i.e. $\boldsymbol{\theta} = \boldsymbol{\theta}_0$) based on the results in Section 2. This allows to give the following lemma.

**Lemma 3.1.** *Under the conditions of 2.1, Corollary 2.1 and Theorem 2.2, the parameters of **Model 1** to **Model 3** are identifiabile for the MLE under Condition **(C7)** and for the GMM under Condition **(C8)**.*

Condition **(C7)** is stronger to assume compared to Condition **(C8)** since, with a few exceptions, the former would the require $f(\boldsymbol{z}|\boldsymbol{\theta})$ to be the density of a Gaussian distribution. On the other hand, Condition **(C8)** is not necessarily too strong since GMM estimators for dependent processes are usually based on the covariance function $\boldsymbol{\varphi}(\boldsymbol{\theta})$, assuming the choice of a positive semi-definite matrix $\boldsymbol{W}$. These conditions are necessary to ensure that one of the basic (and typically difficult to derive) conditions for consistency of an extremum estimator holds: identifiability. Under Conditions **(C7)** and **(C8)**, where the latter is easily verifiable, Section 2 has shown that the covariance function is injective for a wide class of second-order stationary latent time series and spatial models which therefore removes this challenging condition for them. Based on these conditions, the following corollary states the consistency of the MLE and GMM estimators which, for simplicity, we both denote as $\tilde{\boldsymbol{\theta}}$.

**Corollary 3.1.** *Under the conditions of Lemma 3.1 and conditions **(C4)** to **(C6)**, we have that*

$$\tilde{\boldsymbol{\theta}} \xrightarrow{\mathcal{P}} \boldsymbol{\theta}.$$

Having stated the consistency of two popular extremum estimators (i.e. MLE and GMM), the next section discusses the properties of the GMWM for latent models which include the non-stationary processes (T3) and (T4).

### 3.2. GMWM Estimators

The GMWM is a minimum-distance extremum estimator where the reduced form (or auxiliary) parameters are represented by the Wavelet Variance (WV) which represents weighted averages of the Spectral Density Function (SDF) over octave bands. More specifically, the GMWM is defined as follows

$$\hat{\boldsymbol{\theta}} = \underset{\boldsymbol{\theta} \in \boldsymbol{\Theta}}{\operatorname{argmin}} (\hat{\boldsymbol{\nu}} - \boldsymbol{\nu}(\boldsymbol{\theta}))^T \boldsymbol{\Omega} (\hat{\boldsymbol{\nu}} - \boldsymbol{\nu}(\boldsymbol{\theta})) \qquad (3.1)$$

where $\hat{\boldsymbol{\nu}} = [\hat{\nu}_j^2]_{j=1,\ldots,J}$ represents the vector of estimated Haar WV for the $J$ scales of wavelet decomposition (see Percival, 1995), $\boldsymbol{\nu}(\boldsymbol{\theta}) = [\nu_j^2(\boldsymbol{\theta})]_{j=1,\ldots,J}$ represents the vector of WV implied by the model of interest and $\boldsymbol{\Omega}$ is a positive definite weighting matrix chosen in a suitable manner (see Guerrier et al., 2013).

As can be observed from the form of the objective function of the GMWM, the latter depends on $\boldsymbol{\theta}$ uniquely through the theoretical form of the WV $\boldsymbol{\nu}(\boldsymbol{\theta})$ (i.e. the WV implied by the model $F_{\boldsymbol{\theta}}$). Based again on Newey and McFadden (1994), we have identifiability of the parameters $\boldsymbol{\theta}$ through the GMWM if:

1. $\boldsymbol{\nu}(\boldsymbol{\theta}) = \boldsymbol{\nu}(\boldsymbol{\theta}_0)$ if and only if $\boldsymbol{\theta} = \boldsymbol{\theta}_0$; and
2. $\boldsymbol{\Omega}$ is positive definite.

Therefore, we would need to focus on whether the WV $\boldsymbol{\nu}(\boldsymbol{\theta})$ is injective. Let us first focus on the latent models made by the sum of stationary processes. Using the same logic as in Section 2, we would need to prove that the Jacobian matrix $\mathbf{A}(\boldsymbol{\theta}) = \partial/\partial\boldsymbol{\theta}\, \boldsymbol{\nu}(\boldsymbol{\theta})$ is full column rank, keeping in mind that Condition **(C2)** should hold also in this case. However, this direct approach appears difficult to prove when considering the sum of (T6) processes and an approach by steps is needed. This approach discusses the unique mapping of the parameters $\boldsymbol{\theta}$ to the SDF $S_{\boldsymbol{\theta}(f)}$ and then from there to the WV $\boldsymbol{\nu}(\boldsymbol{\theta})$. However, we are forced to assume this last step and therefore we give the following condition.

**(C9)** There is a unique mapping between the SDF $S_{\boldsymbol{\theta}}(f)$ and the WV $\boldsymbol{\nu}(\boldsymbol{\theta})$.

This condition was discussed in the continuous case for the Allan Variance (AV) by Greenhall (1998). The AV has a direct relationship with the Haar WV (i.e. AV $\equiv$ 2WV) and can therefore be considered similar to the WV in that it is derived from another kind of averaging of the SDF. In Greenhall (1998) some very specific cases are discussed in which this mapping is not necessarily unique but, as the author himself claims, these are particular cases that hardly exist in reality. More specifically, if we define the signed SDF as $\Phi(f) = S_{\boldsymbol{\theta}_0}(f) - S_{\boldsymbol{\theta}_1}(f)$, a condition that the author gives for Condition **(C9)** not to be satisfied is given below.

**(C10)** The signed SDF $\Phi(f)$ satisfies $\Phi(2f) = \frac{1}{2}\Phi(f)$.

The reason for this condition resides in the definition of the Haar WV which is adapted from Greenhall (1998) as follows

$$\nu_j^2(\boldsymbol{\theta}) = \frac{1}{2}\int_0^\infty \frac{\sin^4(\tau\pi f)}{(\tau\pi f)^2} S_{\boldsymbol{\theta}}(f) df,$$

where $\tau = 2^j$. As a consequence of this definition, there would not be a unique mapping from the SDF to the WV if $\Delta \equiv \nu_j^2(\boldsymbol{\theta}_0) - \nu_j^2(\boldsymbol{\theta}_1) = 0$ for some $\boldsymbol{\theta}_0 \neq \boldsymbol{\theta}_1$.

Following the proof in Greenhall (1998), using trigonometric inequalities and discarding constants, this would deliver

$$\Delta = \int_0^\infty \frac{\sin^4(\tau\pi f)}{(\tau\pi f)^2}\Phi(f)df = \lim_{n\to-\infty}\int_{2^n}^\infty \frac{\sin^2(\tau\pi f) - \frac{1}{4}\sin^2(2\tau\pi f)}{(\tau\pi f)^2}\Phi(f)df.$$

Rewriting the above as a difference of integrals and by change of variable $u = 2f$ in the second term we have

$$\Delta = \lim_{n\to-\infty}\left[\int_{2^n}^\infty \frac{\sin^2(\tau\pi f)}{(\tau\pi f)^2}\Phi(f)df - \frac{1}{2}\int_{2^{n+1}}^\infty \frac{\sin^2(\tau\pi u)}{(\tau\pi u)^2}\Phi(f)du\right]$$

where, by finally using Condition **(C10)**, we obtain

$$\Delta = \lim_{n\to-\infty}\int_{2^n}^{2^{n+1}} \frac{\sin^2(\tau\pi f)}{(\tau\pi f)^2}\Phi(f)df = 0.$$

This implies that, in the continuous case, there is not a unique mapping from the SDF to the WV when Condition **(C10)** is satisfied. The following corollary gives a result regarding this condition.

**Corollary 3.2.** *The signed SDF of Model 1 does not satisfy Condition (C10).*

The proof of this corollary can be found in Appendix A.4. However, even in the discrete case, Condition **(C10)** gives us the "if" but not the "only if" statement which would allow identifiability. The following conjecture refers to this problem.

**Conjecture 3.1.** *There is always a unique mapping between the SDF and the WV if Condition (C10) is not satisfied.*

This conjecture, whose proof is left for future research, claims that Condition **(C10)** is therefore the only condition which must be verified for there to be a unique mapping from the SDF to the WV and consequently concludes the discussion on Condition **(C9)**. The latter therefore does not appear to be a strong condition to assume also considering additional arguments given further on. Having discussed this condition, we now provide the following corollary.

**Corollary 3.3.** *Under the conditions of Theorem 2.1 and Corollary 2.1, and under Condition (C9), the WV $\nu(\boldsymbol{\theta})$ of Model 1 and Model 2 is injective.*

The proof of this corollary is a direct consequence of Theorem 2.1 by using the argument that a composition of injective functions is itself an injective function. Indeed, Lemmas A.1 and A.2 in Appendix A.5 show that there is a unique mapping between the covariance function and the SDF of these latent models and, assuming Conjecture 3.1 is true, we have that there is a unique mapping between the parameter $\boldsymbol{\theta}$ and the WV $\nu(\boldsymbol{\theta})$. Hence, Condition **(C9)** allows us to avoid proving the full column rank of the Jacobian matrix $\mathbf{A}(\boldsymbol{\theta})$. Nevertheless, this proof is possible for some classes of latent models given below.

**Model 4** $W_t = \sum_{i=6}^{G} X_t^{(i)}$.
**Model 5** $W_t = \sum_{i=1}^{4} X_t^{(i)}$.
**Model 6** $W_t = \sum_{i=3}^{5} X_t^{(i)}$.

**Model 4** is a stationary latent model since it consists in a sum of $K$ (T6) processes while **Model 5** and **Model 6** are classes of latent models that combine stationary and non-stationary processes. Let us start from **Model 4** for which the determinant of the Jacobian $\mathbf{A}(\boldsymbol{\theta})$ is always positive for $G = 10$ (i.e. the sum of four (T6) processes) thereby implying that the WV $\boldsymbol{\nu}(\boldsymbol{\theta})$ is injective for this model. Based on this finding, we give the following conjecture.

**Conjecture 3.2.** *Under the conditions of Theorem 2.1, the determinant of the Jacobian* $\mathbf{A} = \partial/\partial\boldsymbol{\theta}\,\boldsymbol{\nu}(\boldsymbol{\theta})$ *for* **Model 4** *($G < \infty$) is given by*

$$|\mathbf{A}| = \frac{\prod_{i=1}^{K} \upsilon_i^2 \prod_{i<j}^{K} (\rho_i - \rho_j)^4}{\prod_{i=1}^{K} (\rho_i^2 - 1)^2} > 0.$$

As mentioned earlier, this conjecture is verified for a sum of up to four (T6) processes, implying that the Jacobian is of full column rank. Unfortunately, proof by induction is quite hard in this case and is therefore not considered. However this conjecture, along with the results in Greenhall (1998), gives additional support to the argument stating that Condition **(C9)** is not a strong one to assume. Finally, using the same approach based on the Jacobian matrix $\mathbf{A}(\boldsymbol{\theta})$, let us give the following lemma.

**Lemma 3.2.** *Under Conditions (C1) and (C2), the WV $\boldsymbol{\nu}(\boldsymbol{\theta})$ of* **Model 5** *and* **Model 6** *is injective.*

This lemma, whose proof is given in Appendix A.6, allows to combine some stationary processes with the non-stationary processes (T3) and (T4) considered in this paper. Indeed, excluding the (T6) process, the parameters of each stationary process can be identified together with those of the non-stationary ones. With these results, we can finally give the following lemma.

**Lemma 3.3.** *The GMWM is able to identify the parameters of*

- **Model 1** *and* **Model 2** *under the conditions of Corollary 3.3.*
- **Model 4** *under Condition (C9).*
- **Model 5** *and* **Model 6** *under the conditions of Lemma 3.2.*

This corollary therefore summarizes the results of this section. Indeed, the parameters of **Model 1** and **Model 2** are identifiable for the GMWM if we assume Condition **(C9)** which is claimed to be true for these latent models in Conjecture 3.1. On the other hand, Conjecture 3.2 claims that the parameters of **Model 4** are also identifiable for the GMWM while those of **Model 5** and **Model 6** are generally identifiable. Considering these results, the consistency of the GMWM is given in the following lemma.

**Corollary 3.4.** *For the processes and conditions considered in Corollary 3.3, and assuming Conditions (C4) to (C6), we have that*

$$\hat{\boldsymbol{\theta}} \xrightarrow{\mathcal{P}} \boldsymbol{\theta}.$$

The GMWM is therefore consistent for different classes of latent models based on the above results. The case of a latent models including a sum of (T6) processes with processes (T3) and (T4) has not been investigated here since the approach used to prove Corollary 3.3 is based on the SDF which is not defined for non-stationary processes. However, if processes (T1) to (T5) are only included once in a latent model, the results of this section strongly suggest that the parameters of any latent model made by the combination of the time series models considered in this paper, excluding a combination of (T1) and/or (T2) with (T5), is identifiable through the WV $\boldsymbol{\nu}(\boldsymbol{\theta})$. Corollary 3.3 and Lemma 3.2 support this idea based on the intuitive argument that the WV of the non-stationary processes (T3) and (T4) increases steadily at the larger scales which cannot be in any way approximated by the stationary processes considered here since their WV decreases at these scales.

The use of the GMWM for the estimation of the parameters of latent spatial models would undergo the same type of conditions and arguments given above for latent time series models. Indeed, a similar statement to that of Corollary 3.3 could be envisaged for **Model 3** consisting in the sum of either spatial models (S1) or (S2). However, it is not clear if the same arguments for the mapping of the covariance function to the SDF can be used and directly proving the full-column rank of the Jacobian matrix $\mathbf{A}(\boldsymbol{\theta}) = \partial/\partial\boldsymbol{\theta}\,\boldsymbol{\nu}(\boldsymbol{\theta})$ is equally as challenging as for the latent time series models made by a sum of (T6) processes.

### 3.3. Numerical Example

To conclude this section, we provide a brief numerical example in which we compare the performance of the GMM and GMWM estimators when estimating the parameters of two latent models on different sample sizes. These two models are defined as follows:

**Latent Model 1** $W_t = X_t^{(1)} + X_t^{(2)} + \sum_{i=6}^{7} X_t^{(i)}$ with parameter vector $[\rho_1 \; \upsilon_1^2 \; \rho_2 \; \upsilon_2^2 \; \sigma^2 \; Q^2]^T = [0.99 \; 1 \; 0.85 \; 15 \; 2 \; 4]$.

**Latent Model 2** $W_t = X_t^{(1)} + X_t^{(2)} + X_t^{(5)} + X_t^{(6)}$ with parameter vector $[\rho_1 \; \upsilon_1^2 \; \varrho \; \varsigma^2 \; \sigma^2 \; Q^2]^T = [0.99 \; 1 \; 0.85 \; 15 \; 2 \; 4]$.

**Latent Model 1** is made by the sum of two (T6) processes with a (T1) and a (T2) process while **Latent Model 2** is the same except that the second autoregressive process is replaced by a (T5) process. From the discussions in Section 2, **Latent Model 2** could have issues of identifiability since process (T5) is included in a model with processes (T1) and (T2).

These models were estimated 500 times for each sample size going from 100 to 1 million and the Mean Squared Error (MSE) of the estimators was computed for each parameter. The results are given in Figure 1 where the logarithmic transform of the MSE for each parameter is given along with its confidence intervals (obtained via bootstrap on the estimated parameters) for growing sample sizes.

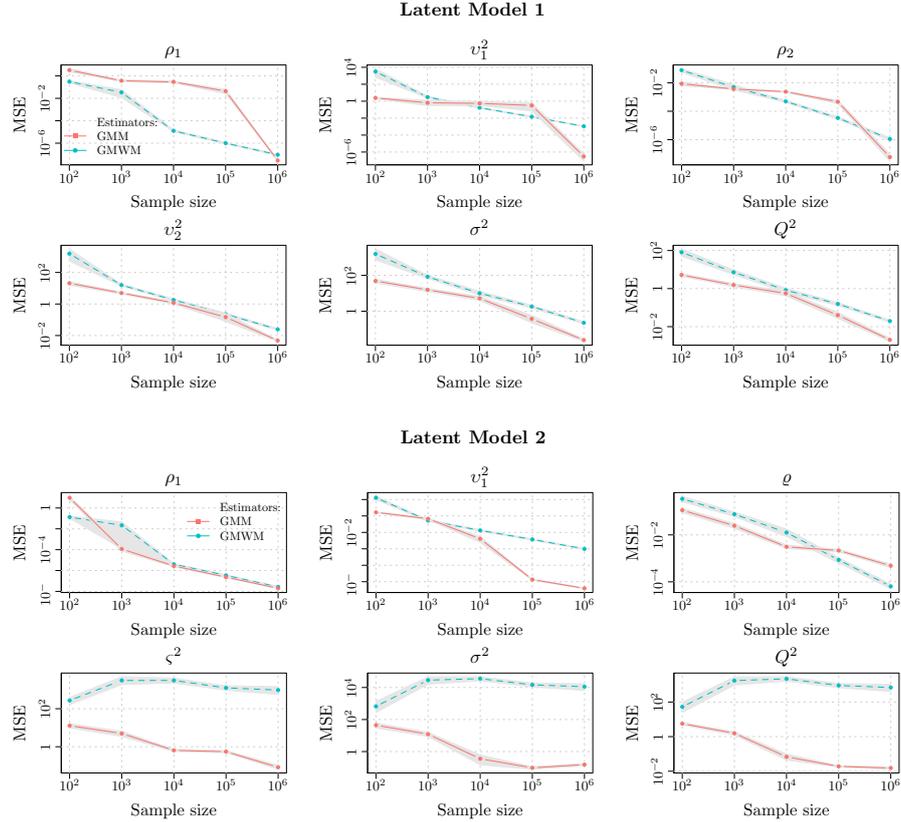

FIG 1. *Logarithm of the MSE for the GMM (red line) and GMWM (green line) on Latent Model 1 (top part of figure) and Latent Model 2 (bottom part of figure). The grey shaded areas represent the confidence intervals of the MSE for each estimator.*

For consistent estimators we would expect the logarithmic transform of the MSE to decrease monotonically as the sample size increases, which is roughly the case for the parameters of **Latent Model 1**. However, this is not the case for the MSE of **Latent Model 2** where it does not have a clear monotonic decrease for the last three parameters that belong to the (T1), (T2) and (T5) processes. Indeed, in this case the GMWM seems to remain roughly constant as the sample size increases while the GMM appears to reach a lower bound for the parameters of processes (T1) and (T2). These results seem therefore to support some of the findings highlighted in the previous sections.

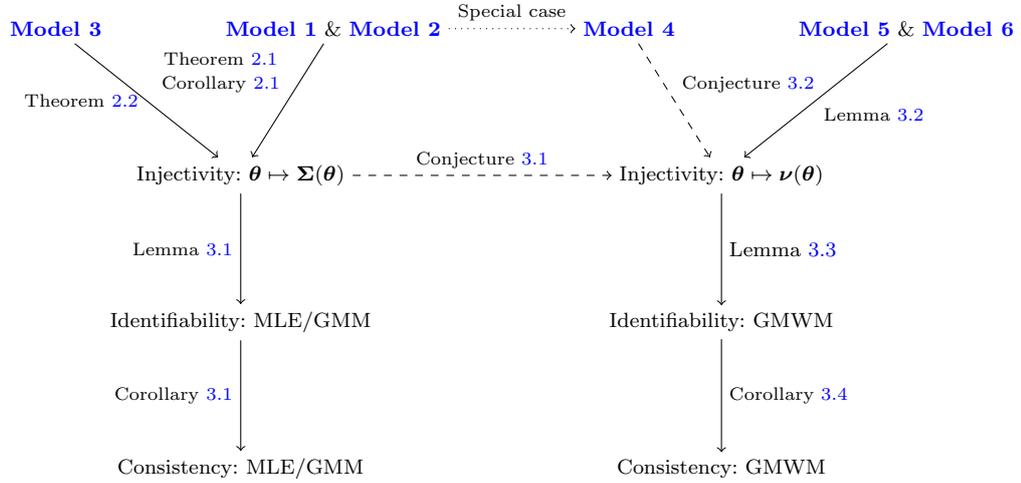

FIG 2. *Summary of the results on identifiability and consistency for the MLE, GMM and GMWM estimators concerning the classes of latent models considered in this paper.*

## 4. Conclusion

The results on the identifiability of latent models presented in this paper are summarized in Figure 2 and considerably widen the spectrum of identifiable models in dependent data settings beyond standard ARMA models. As a result, the widespread use of these types of models is supported by the findings of this paper which also highlighted the latent models which are not identifiable, thereby avoiding an inappropriate estimation procedure. Moreover, we extended the idea of latent models to a group of spatial models for which we also proved that a combination thereof is identifiable, delivering a flexible class of models for spatial data. Finally, we showed how these results considerably reduce the conditions needed for consistency for some extremum estimators, such as the MLE and GMM, while the addition of non-stationary models in this framework and their identifiability was discussed for the GMWM estimators.

**Appendix A: Proofs**

*A.1. Proof of Theorem 2.1*

Let us consider **Model 1** and, to start, let us just consider a sum of independent (T6) processes, which we denote simply as $Y_t = \sum_{i=1}^{K} X_t^{(i)}$, for which a part of the ACVF sequence $(\varphi_{\boldsymbol{\theta}}(h))_{h=0}^{\infty}$ with positive lags is given by

$$\varphi_{\boldsymbol{\theta}}(h) = \sum_{i=1}^{K} \rho_i^h \frac{v_i^2}{1-\rho_i^2},$$

where $\boldsymbol{\theta} = [\rho_1 \ v_1^2 \cdots \rho_K \ v_K^2]$ represents the parameter vector containing the parameters of the $K$ processes. The derivatives with respect to $\rho_i$ and $v_i^2$ (i.e. the parameters of the $i^{th}$ (T6) process) are respectively

$$\begin{aligned}
\gamma_1^i(h) \equiv \frac{\partial}{\partial \rho_i} \varphi_{\boldsymbol{\theta}}(h) &= \frac{h\rho_i^{h-1}(1-\rho_i^2) + 2\rho_i^{h+1}}{(1-\rho_i^2)^2} v_i^2 \\
&= \frac{\rho_i^h(h + (2-h)\rho_i^2)}{\rho_i(1-\rho_i^2)^2} v_i^2 \\
\gamma_2^i(h) \equiv \frac{\partial}{\partial v_i^2} \varphi_{\boldsymbol{\theta}}(h) &= \frac{\rho_i^h}{(1-\rho_i^2)}
\end{aligned}$$

which exist based on the parameter values defined for process (T6). These thereby deliver a Jacobian matrix A whose first $2K$ rows are

$$A = \begin{pmatrix} \gamma_1^1(0) & \gamma_2^1(0) & \cdots & \gamma_1^K(0) & \gamma_2^K(0) \\ \gamma_1^1(1) & \gamma_2^1(1) & \cdots & \gamma_1^K(1) & \gamma_2^K(1) \\ \vdots & \vdots & \vdots & \vdots & \vdots \\ \gamma_1^1(H) & \gamma_2^1(H) & \cdots & \gamma_1^K(H) & \gamma_2^K(H) \end{pmatrix}$$

where $H = 2K - 1$. To simplify notation, let us define $\boldsymbol{\gamma}_1^i(h) = [\gamma_1^i(h)]_{h=0,\ldots,H}$ and $\boldsymbol{\gamma}_2^i(h) = [\gamma_2^i(h)]_{h=0,\ldots,H}$ as being the columns of the matrix A and $|A|$ as being the determinant of A. Taking the determinant of the matrix A, we perform some column permutations and operations. First, for more clarity in the proof, we permute the columns to obtain the following matrix determinant

$$|A| = (-1)^{K-1} |\boldsymbol{\gamma}_1^1(h) \ \cdots \ \boldsymbol{\gamma}_1^K(h) \ \boldsymbol{\gamma}_2^1(h) \ \cdots \ \boldsymbol{\gamma}_2^K(h)|.$$

Next we multiply each column by a different constant which leaves us with the modified columns

$$\begin{aligned}
\tilde{\boldsymbol{\gamma}}_1^i(h) &\equiv \boldsymbol{\gamma}_1^i(h) \frac{\rho_i(1-\rho_i^2)^2}{v_i^2} = \rho_i^h(h + (2-h)\rho_i^2) = \rho_i^h h + (2-h)\rho_i^{h+2} \\
\tilde{\boldsymbol{\gamma}}_2^i(h) &\equiv \boldsymbol{\gamma}_2^i(h)(1-\rho_i^2) = \rho_i^h
\end{aligned}$$

thereby leaving us with the determinant

$$|A| = \underbrace{(-1)^{K-1} \prod_{i=1}^{K} \left( \frac{v_i^2}{\rho_i(1-\rho_i^2)^3} \right)}_{c} |\tilde{\gamma}_1^1(h) \;\cdots\; \tilde{\gamma}_1^K(h) \; \tilde{\gamma}_2^1(h) \;\cdots\; \tilde{\gamma}_2^K(h)|$$

where $c \neq 0$. Now let us express $\tilde{\gamma}_1^i(h)$ as $\boldsymbol{\delta}_1^i(h) + \boldsymbol{\delta}_2^i(h)$ where $\boldsymbol{\delta}_1^i(h) = [\rho_i^h h]_{h=0,\ldots,H}$ and $\boldsymbol{\delta}_2^i(h) = [(2-h)\rho_i^{h+2}]_{h=0,\ldots,H}$. Furthermore, let us define the set $\mathcal{S} = \{\{a_1, \ldots, a_K\} | a_i \in \{1,2\}, i = 1, \ldots, K\}$ with $s \in \mathcal{S}$ being an element of this set and $s(i) = a_i$ being the $i^{th}$ element of $s$ (note that the cardinality of $\mathcal{S}$, and therefore $s$, is $K$). When using the "determinant as a sum of determinants" rule, we can split each column $\tilde{\gamma}_1^i(h)$ which, starting with $\tilde{\gamma}_1^1(h)$, creates a sum of the final nodes of a binary tree as follows

$$\begin{aligned}
|A| &= c(|\boldsymbol{\delta}_1^1(h) \; \tilde{\gamma}_1^2(h) \;\cdots\; \tilde{\gamma}_1^K(h) \; \tilde{\gamma}_2^1(h) \;\cdots\; \tilde{\gamma}_2^K(h)| \\
&\quad + |\boldsymbol{\delta}_2^1(h) \; \tilde{\gamma}_1^2(h) \;\cdots\; \tilde{\gamma}_1^K(h) \; \tilde{\gamma}_2^1(h) \;\cdots\; \tilde{\gamma}_2^K(h)|) \\
&= c(|\boldsymbol{\delta}_1^1(h) \; \boldsymbol{\delta}_1^2(h) \; \tilde{\gamma}_1^3(h) \;\cdots\; \tilde{\gamma}_1^K(h) \; \tilde{\gamma}_2^1(h) \;\cdots\; \tilde{\gamma}_2^K(h)| \\
&\quad + |\boldsymbol{\delta}_1^1(h) \; \boldsymbol{\delta}_2^2(h) \; \tilde{\gamma}_1^3(h) \;\cdots\; \tilde{\gamma}_1^K(h) \; \tilde{\gamma}_2^1(h) \;\cdots\; \tilde{\gamma}_2^K(h)| \\
&\quad + |\boldsymbol{\delta}_2^1(h) \; \boldsymbol{\delta}_1^2(h) \; \tilde{\gamma}_1^3(h) \;\cdots\; \tilde{\gamma}_1^K(h) \; \tilde{\gamma}_2^1(h) \;\cdots\; \tilde{\gamma}_2^K(h)| \\
&\quad + |\boldsymbol{\delta}_2^1(h) \; \boldsymbol{\delta}_2^2(h) \; \tilde{\gamma}_1^3(h) \;\cdots\; \tilde{\gamma}_1^K(h) \; \tilde{\gamma}_2^1(h) \;\cdots\; \tilde{\gamma}_2^K(h)|) \\
&= \cdots \\
&= c \sum_{s \in \mathcal{S}} |\boldsymbol{\delta}_{s(1)}^1(h) \cdots \boldsymbol{\delta}_{s(K)}^K(h) \; \tilde{\gamma}_2^1(h) \;\cdots\; \tilde{\gamma}_2^K(h)|.
\end{aligned}$$

Notice that each column $\boldsymbol{\delta}_{s(i)}^i(h)$ with $s(i) = 2$ (i.e. $\boldsymbol{\delta}_2^i(h)$) can be re-expressed as $\boldsymbol{\delta}_2^i(h) = [(2-h)\rho_i^h \rho_i^2]_{h=0,\ldots,H}$ and therefore, denoting $\tilde{\boldsymbol{\delta}}_2^i(h) = [(2-h)\rho_i^h]_{h=0,\ldots,H}$, we have that the determinant becomes

$$|A| = c \sum_{s \in \mathcal{S}} \underbrace{\left( \prod_{i=1}^{K} (\mathbb{1}_{s(i)=2}\rho_i^2 + \mathbb{1}_{s(i)=1}) \right)}_{c_s} |\tilde{\boldsymbol{\delta}}_{s(1)}^1(h) \cdots \tilde{\boldsymbol{\delta}}_{s(K)}^K(h) \; \tilde{\gamma}_2^1(h) \;\cdots\; \tilde{\gamma}_2^K(h)|.$$

where $\tilde{\boldsymbol{\delta}}_{s(i)}^i(h) = \mathbb{1}_{s(i)=2}\tilde{\boldsymbol{\delta}}_2^i(h) + \mathbb{1}_{s(i)=1}\boldsymbol{\delta}_1^i(h)$ and $\mathbb{1}_{s(i)=z}$ represents the indicator function which takes the value 1 if $s(i) = z$ and 0 otherwise. Following the same procedure as before, let us express $\tilde{\boldsymbol{\delta}}_2^i(h) = \boldsymbol{\delta}_3^i(h) + \boldsymbol{\delta}_4^i(h)$, where $\boldsymbol{\delta}_3^i(h) = [2\rho_i^h]_{h=0,\ldots,H} = 2\tilde{\gamma}_2^i(h)$ and $\boldsymbol{\delta}_4^i(h) = [-\rho_i^h h]_{h=0,\ldots,H} = -\boldsymbol{\delta}_1^i(h)$. If we split the columns $\tilde{\boldsymbol{\delta}}_2^i(h)$ into the sub-determinants (in the same way as for $\tilde{\gamma}_1^i(h)$) we have a binary tree for which each node has two children of which the one which includes the column $\boldsymbol{\delta}_3^i(h)$ will be null since this column is a linear function of the column $\tilde{\gamma}_2^i(h)$. This implies that the only term which remains in each split is $\boldsymbol{\delta}_4^i(h)$ which is simply the negative of $\boldsymbol{\delta}_1^i(h)$ and there will be as many $\boldsymbol{\delta}_4^i(h)$ columns as the original $\tilde{\boldsymbol{\delta}}_2^i(h)$ columns. Given this structure, we have that the

determinant becomes

$$|\mathrm{A}| = c \underbrace{\sum_{s \in \mathcal{S}} c_s (-1)^{\sum_{i=1}^{K} \mathbb{1}_{s(i)=2}}}_{c_{s*}} \underbrace{|\boldsymbol{\delta}_1^1(h) \cdots \boldsymbol{\delta}_1^K(h) \ \tilde{\boldsymbol{\gamma}}_2^1(h) \ \cdots \ \tilde{\boldsymbol{\gamma}}_2^K(h)|}_{\lambda}.$$

We therefore have that the determinant $\lambda$ is the determinant of the following matrix

$$B = \begin{pmatrix} \rho_1^0 0 & \rho_2^0 0 & \cdots & \rho_K^0 0 & \rho_1^0 & \rho_2^0 & \cdots & \rho_K^0 \\ \rho_1^1 1 & \rho_2^1 1 & \cdots & \rho_K^1 1 & \rho_1^1 & \rho_2^1 & \cdots & \rho_K^1 \\ \vdots & \vdots & \vdots & \vdots & \vdots & \vdots & \vdots & \vdots \\ \rho_1^H H & \rho_2^H H & \cdots & \rho_K^H H & \rho_1^H & \rho_2^H & \cdots & \rho_K^H \end{pmatrix}$$

which is clearly full column rank. Indeed, the last $K$ columns are those of a Vandermonde matrix with distinct elements $\rho_i \neq \rho_j$, $\forall i \neq j$, implying that the columns are linearly independent, while the first $K$ columns are a Vandermonde matrix whose rows are multiplied by distinct constants thereby implying that also these columns are linearly independent. Moreover, the first $K$ columns can be seen as functions of the last $K$ columns which cannot however be expressed as a linear combination of the others. Therefore there is no column that can be expressed as a linear combination of the others implying that $\lambda \neq 0$. As for the term $c_{s*}$ it can be interpreted in a geometric manner since it represents the $K$-dimensional volume of a hyperrectangle with sides $(1 - \rho_i^2)$ given by

$$c_{s*} = \prod_{i=1}^{K} (1 - \rho_i^2),$$

where $0 < \prod_{i=1}^{K}(1-\rho_i^2) < 1$. If we permute the columns so as to respect the order of the parameters in the original vector $\boldsymbol{\theta}$, this finally delivers the determinant

$$|\mathrm{A}| = \prod_{i=1}^{K} \left( \frac{v_i^2}{\rho_i(1-\rho_i^2)^2} \right) \lambda \neq 0,$$

thereby implying that the Jacobian matrix A is of full column rank and, consequently, that the ACVF $\varphi_{\boldsymbol{\theta}}(h)$ of a sum of $K$ (T6) processes is injective. Notice that if we were to add processes (T1) and (T2) to the above procedure, the only detail that would change is the matrix B. Indeed, considering the derivatives of their respective ACVFs, if they were inserted in the original matrix A these columns would remain unaffected by the operations thereby delivering the matrix

$$\tilde{B} = \begin{pmatrix} \rho_1^0 0 & \rho_2^0 0 & \cdots & \rho_K^0 0 & \rho_1^0 & \rho_2^0 & \cdots & \rho_K^0 & 1 & 2 \\ \rho_1^1 1 & \rho_2^1 1 & \cdots & \rho_K^1 1 & \rho_1^1 & \rho_2^1 & \cdots & \rho_K^1 & 0 & -1 \\ \vdots & \vdots & \vdots & \vdots & \vdots & \vdots & \vdots & \vdots & \vdots & 0 \\ \vdots & \vdots & \vdots & \vdots & \vdots & \vdots & \vdots & \vdots & \vdots & \vdots \\ \rho_1^H H & \rho_2^H H & \cdots & \rho_K^H H & \rho_1^H & \rho_2^H & \cdots & \rho_K^H & 0 & 0 \end{pmatrix}$$

which is also clearly full rank with determinant $\tilde{\lambda} \neq 0$. This implies that the ACVF of a latent process made by the sum of $K$ (T6) processes, a (T1) process and a (T2) process is also injective through the parameter vector $\boldsymbol{\theta}$. $\square$

### A.2. Proof of injection through an ARMA(2,1)

Let us denote the parameters of an ARMA(2,1) process as $\tilde{\boldsymbol{\theta}} = [\tilde{\rho}_1 \ \tilde{\rho}_2 \ \tilde{\varrho} \ \tilde{v}^2]^T$. Given the reparametrization of the sum of two (T6) processes to an ARMA(2,1) process given in Hamilton (1994) (equation 4.7.26, Section 4.7), we have that the parameters of an ARMA(2,1) process are given by

$$\tilde{\boldsymbol{\theta}} = \begin{bmatrix} \tilde{\rho}_1 \\ \tilde{\rho}_2 \\ \tilde{\varrho} \\ \tilde{v}^2 \end{bmatrix} = \begin{bmatrix} \rho_1 + \rho_2 \\ -\rho_1 \rho_2 \\ \frac{\rho_1 v_1^2 + \rho_2 v_2^2}{v_1^2 + v_2^2} \\ v_1^2 + v_2^2 \end{bmatrix} = \boldsymbol{g}(\boldsymbol{\theta})$$

where $\boldsymbol{\theta} = [\rho_1 \ v_1^2 \ \rho_2 \ v_2^2]^T$ is the vector of parameters for the sum of two (T6) processes. By taking the Jacobian $A = \partial/\partial \boldsymbol{\theta} \ \boldsymbol{g}(\boldsymbol{\theta})$ and using the reasoning given in (2.1), we have that the determinant of this matrix is given by

$$|A| = \frac{(\rho_1 - \rho_2)^2}{v_1^2 + v_2^2}$$

which is always positive thereby implying that the matrix A is of full rank and that there is a unique mapping from the parameters of a sum of two (T6) processes to those of an ARMA(2,1) process. $\square$

### A.3. Proof of Theorem 2.2

This proof follows the procedure carried out to prove Theorem 2.1. Having said this, let us consider **Model 3** and let

$$\gamma_{\boldsymbol{\theta}_i}^c(d) = \sigma_i^2 \exp\left(-\left(\frac{d}{\phi_i}\right)^c\right) \quad (A.1)$$

be the covariance function of the $i^{th}$ spatial model with parameters $\phi_i$ and $\sigma_i^2$ and where $d$ is a distance between two point $k$ and $k'$, while $c \in \mathbb{N}^+$ is a known constant. Special cases of this spatial model are the Exponential ($c = 1$) and Gaussian ($c = 2$) models. Defining $\omega_i = \exp(-1/\phi_i^c)$, the derivatives with respect to $\phi_i$ and $\sigma_i^2$ are respectively

$$\alpha_i(d) \equiv \frac{\partial}{\partial \phi_i^2} \gamma_{\boldsymbol{\theta}_i}^c(d) \quad = \quad \sigma_i^2 \omega_i^{d^c} c \frac{d^{c-1}}{\phi^{c-1}} \frac{d}{\phi^2} = \underbrace{\frac{\sigma_i^2 c}{\phi_i^{c+1}}}_{\Delta_i^c} \omega_i^{d^c} d^c$$

$$\beta_i(d) \equiv \frac{\partial}{\partial \sigma_i^2} \gamma^c_{\boldsymbol{\theta}_i}(d) = \omega_i^{d^c}.$$

With $\boldsymbol{\theta} = [\theta_1 \cdots \theta_K]^T = [\phi_1 \ \sigma_1^2 \ \cdots \ \phi_K \ \sigma_K^2]^T$, these derivatives thereby deliver a Jacobian matrix A whose first $2K$ rows are

$$A = \begin{pmatrix} \alpha_1(d_0) & \beta_1(d_0) & \cdots & \alpha_K(d_0) & \beta_K(d_0) \\ \alpha_1(d_1) & \beta_1(d_1) & \cdots & \alpha_K(d_1) & \beta_K(d_1) \\ \vdots & \vdots & \vdots & \vdots & \vdots \\ \alpha_1(d_H) & \beta_1(d_H) & \cdots & \alpha_K(d_H) & \beta_K(d_H) \end{pmatrix}$$

where $H = 2K - 1$ and $d_i < d_j$, $\forall i < j$ with $d_0 = 0$. Let us define the columns of A as $\boldsymbol{\alpha}_i(d) = [\alpha_i(d)]_{d=d_0,\ldots,d_H}$ and $\boldsymbol{\beta}_i(d) = [\beta_i(d)]_{d=d_0,\ldots,d_H}$ respectively such that the determinant of matrix A is denoted as

$$A = |\boldsymbol{\alpha}_1(d) \ \boldsymbol{\beta}_1(d) \ \cdots \ \boldsymbol{\alpha}_K(d) \ \boldsymbol{\beta}_K(d)|.$$

Rearranging the columns we obtain

$$|A| = (-1)^K |\boldsymbol{\alpha}_1(d) \ \cdots \ \boldsymbol{\alpha}_K(d) \ \boldsymbol{\beta}_1(d) \ \cdots \ \boldsymbol{\beta}_K(d)|.$$

We know that each column $\boldsymbol{\alpha}_i(d)$ is multiplied by a constant $\Delta_i^c$ such that the matrix determinant becomes

$$|A| = \underbrace{(-1)^K \prod_{i=1}^K \Delta_i^c}_{\neq 0} \underbrace{\begin{vmatrix} \omega_1^{d_0^c} d_0^c & \cdots & \omega_K^{d_0^c} d_0^c & \omega_1^{d_0^c} & \cdots & \omega_K^{d_0^c} \\ \omega_1^{d_1^c} d_1^c & \cdots & \omega_K^{d_1^c} d_1^c & \omega_1^{d_1^c} & \cdots & \omega_K^{d_1^c} \\ \vdots & \vdots & \vdots & \vdots & \vdots \\ \omega_1^{d_H^c} d_H^c & \cdots & \omega_K^{d_H^c} d_H^c & \omega_1^{d_H^c} & \cdots & \omega_K^{d_H^c} \end{vmatrix}}_{\lambda}$$

where $\lambda$ therefore represents the determinant of the matrix

$$B_c(d_0^H) = \begin{pmatrix} \omega_1^{d_0^c} d_0^c & \cdots & \omega_K^{d_0^c} d_0^c & \omega_1^{d_0^c} & \cdots & \omega_K^{d_0^c} \\ \omega_1^{d_1^c} d_1^c & \cdots & \omega_K^{d_1^c} d_1^c & \omega_1^{d_1^c} & \cdots & \omega_K^{d_1^c} \\ \vdots & \vdots & \vdots & \vdots & \vdots \\ \omega_1^{d_H^c} d_H^c & \cdots & \omega_K^{d_H^c} d_H^c & \omega_1^{d_H^c} & \cdots & \omega_K^{d_H^c} \end{pmatrix}$$

which is clearly full column rank by the same argument used for matrix B in Appendix A.1. Therefore there is no column that can be expressed as a linear combination of the others implying that $\lambda \neq 0$. Based on this, we know that there is a unique mapping of the parameter vector $\boldsymbol{\theta}$ to the covariance function $\gamma^c_{\boldsymbol{\theta}}(d)$.

### A.4. Proof of Corollary 3.2

For a sum of $K$ (T6) processes the SDF is given by

$$S_{\boldsymbol{\theta}}(f) = \sum_{j=1}^K \frac{v_j^2}{1 - 2\rho_j \cos(2\pi f) + \rho_j^2}.$$

If we consider Condition (C10), with $\rho_j$ and $v_j^2$ denoting the elements of the true parameter vector $\boldsymbol{\theta}_0$, we have

$$\sum_{j=1}^{K} \frac{v_j^2}{1 - 2\rho_j \cos(4\pi f) + \rho_j^2} - \sum_{j=1}^{K} \frac{\tilde{v}_j^2}{1 - 2\tilde{\rho}_j \cos(4\pi f) + \tilde{\rho}_j^2}$$
$$= \sum_{j=1}^{K} \frac{v_j^2}{2(1 - 2\rho_j \cos(2\pi f) + \rho_j^2)} - \sum_{j=1}^{K} \frac{\tilde{v}_j^2}{2(1 - 2\tilde{\rho}_j \cos(2\pi f) + \tilde{\rho}_j^2)}$$

where $\tilde{\rho}_j$ and $\tilde{v}_j^2$ denote the elements of a parameter vector $\boldsymbol{\theta}_1 \neq \boldsymbol{\theta}_0$. This condition is clearly not respected since the left side of the equation cannot be expressed as the right side. Indeed, we would need to obtain $2(1 - 2\rho_j \cos(2\pi f) + \rho_j^2)$ from the expression $(1 - 2\rho_j \cos(4\pi f) + \rho_j^2)$ which is not possible. If we add processes (T1) and (T2) to this condition and denote the spectral density of the $i^{th}$ (T6) process as $\alpha_i$ we have

$$\sigma^2 + \frac{4Q^2}{\tau} \sin^2(\pi 2 f \tau) + \sum_{i=1}^{K} \alpha_i - \tilde{\sigma}^2 - \frac{4\tilde{Q}^2}{\tau} \sin^2(\pi 2 f \tau) - \sum_{i=1}^{K} \tilde{\alpha}_i$$
$$= \frac{\sigma^2}{2} + \frac{2Q^2}{\tau} \sin^2(\pi f \tau) + \frac{1}{2} \sum_{i=1}^{K} \alpha_i - \frac{\tilde{\sigma}^2}{2} - \frac{2\tilde{Q}^2}{\tau} \sin^2(\pi f \tau) - \frac{1}{2} \sum_{i=1}^{K} \tilde{\alpha}_i.$$

Also in this case, it is straightforward to see that Condition (C10) is not respected either since the spectral density of process (T1) does not depend on $f$ and is constant while that of process (T2) cannot be re-expressed in this form. □

### A.5. Unique mapping from the covariance function to the spectral density

**Lemma A.1.** *There is a unique mapping between the covariance function $\varphi_{\boldsymbol{\theta}}(h)$ and the spectral density $S_{\boldsymbol{\theta}}(f)$ for a process made by the sum of a (T1) process, a (T2) and K (T6) processes, $\forall K < \infty$.*

For there to be a unique mapping from the ACVF to the spectral density $S_{\boldsymbol{\theta}}(f)$ we need the sequence of autocovariances $(\varphi_{\boldsymbol{\theta}})(h))_{h=-\infty}^{+\infty}$ to be absolutely summable (see, for example, Proposition 6.1 in Hamilton, 1994). Given the parameter values for $\rho_k$, each of the $K$ (T6) processes can be expressed as a causal linear process $Y_t^{(k)} = \sum_{j=0}^{\infty} \psi_j^{(k)} \epsilon_{t-j}^{(k)}$ where $\psi_j^{(k)}$ are fixed coefficients which respect $\sum_{j=0}^{\infty} |\psi_j^{(k)}| < \infty$ and $\epsilon_t^{(k)} \stackrel{iid}{\sim} \mathcal{N}(0, \kappa_k^2)$. Their sum gives

$$\sum_{k=1}^{K} \sum_{j=0}^{\infty} \psi_j^{(k)} \epsilon_{t-j}^{(k)} = \sum_{j=0}^{\infty} \sum_{k=1}^{K} \psi_j^{(k)} \epsilon_{t-j}^{(k)} \tag{A.2}$$

and, given that the processes are independent by condition **(C1)**, we have that

$$\sum_{k=1}^{K} \psi_j^{(k)} \epsilon_{t-j}^{(k)} \sim \mathcal{N}\left(0, \underbrace{\sum_{k=1}^{K} (\psi_j^{(k)})^2 \kappa_k^2}_{\delta^2}\right).$$

We have that $\delta^2$ can be re-expressed as

$$\begin{aligned}
\delta^2 &= (\psi_j^{(1)})^2 \kappa_1^2 + \kappa_1^2 \sum_{k=2}^{K} \frac{\kappa_k^2}{\kappa_1^2} (\psi_j^{(k)})^2 \\
&= \underbrace{\left((\psi_j^{(1)})^2 + \sum_{k=2}^{K} \frac{\kappa_k^2}{\kappa_1^2} (\psi_j^{(k)})^2\right)}_{\tilde{\psi}_j^2} \kappa_1^2.
\end{aligned}$$

We therefore have the equivalence of processes

$$\sum_{k=1}^{K} \psi_j^{(k)} \epsilon_{t-j}^{(k)} \iff \tilde{\psi}_j u_{t-j}, \ \ u_t \sim \mathcal{N}(0, \kappa_1^2)$$

which implies that (A.2) can be written as $\sum_{j=0}^{\infty} \tilde{\psi}_j u_{t-j}$. If we define $\boldsymbol{z} = \left[\psi_j^{(1)} \ \frac{\kappa_2}{\kappa_1} \psi_j^{(2)} \ \frac{\kappa_3}{\kappa_1} \psi_j^{(3)} \cdots \frac{\kappa_K}{\kappa_1} \psi_j^{(K)}\right]$, this gives us

$$\begin{aligned}
\sum_{j=0}^{\infty} |\tilde{\psi}_j| &= \sum_{j=0}^{\infty} \sqrt{(\psi_j^{(1)})^2 + \sum_{k=2}^{K} \frac{\kappa_k^2}{\kappa_1^2} (\psi_j^{(k)})^2} = \sum_{j=0}^{\infty} \|\boldsymbol{z}\|_2 \\
&\leq \sum_{j=0}^{\infty} \|\boldsymbol{z}\|_1 = \sum_{j=0}^{\infty} |\psi_j^{(1)}| + \frac{\kappa_2}{\kappa_1} \sum_{j=0}^{\infty} |\psi_j^{(2)}| + \cdots + \frac{\kappa_K}{\kappa_1} \sum_{j=0}^{\infty} |\psi_j^{(K)}| < \infty
\end{aligned}$$

where $\|\boldsymbol{z}\|_2$ and $\|\boldsymbol{z}\|_1$ represent the $L^2$ norm and $L^1$ norm of $\boldsymbol{z}$ respectively. This result implies that a sum of $K$ (T6) processes have absolutely summable linear coefficients $(\tilde{\psi}_j)_{j=0}^{\infty}$ which implies that its autocovariances $(\varphi_{\boldsymbol{\theta}})(h))_{h=-\infty}^{+\infty}$ are also absolutely summable. If we add the processes (T1) and (T2) to this reasoning, we have that their ACVF are given by

$$\varphi_{\sigma^2}(h) = \begin{cases} \sigma^2 & \text{if } h = 0 \\ 0 & \text{if } h > 0 \end{cases}$$

and

$$\varphi_{Q^2}(h) = \begin{cases} 2Q^2 & \text{if } h = 0 \\ -Q^2 & \text{if } h = 1 \\ 0 & \text{if } h > 1 \end{cases}$$

where $\varphi_{\sigma^2}(h)$ and $\varphi_{Q^2}(h)$ represent the ACVF of the (T1) and (T2) process respectively. If added to the ACVF of the sum of $K$ (T6) processes described earlier, it is straightforward to see that the resulting ACVF sequence is absolutely summable as well, thus concluding the proof. □

**Lemma A.2.** *There is a unique mapping between the covariance function $\varphi_{\boldsymbol{\theta}}(h)$ and the spectral density $S_{\boldsymbol{\theta}}(f)$ for a process made by the sum of a (T5) process and $K$ (T6) processes, $\forall K < \infty$.*

The proof uses the same arguments as the proof of Lemma A.1. □

### *A.6. Proof of Lemma 3.2*

The proof is straightforward for both classes of processes. Considering the process $Y_t = \sum_{i=1}^{4} X_t^{(i)}$, if we take the first four consecutive WV scales (i.e. $\boldsymbol{\nu}(\boldsymbol{\theta_1}) = [\nu_1^2(\boldsymbol{\theta_1}), \ldots, \nu_4^2(\boldsymbol{\theta_1})]$), we have that the determinant of the relative Jacobian matrix is $|A_{4,4}| = {}^{2205\omega}\!/\!{}_{4096}$ which implies that it is of full rank. Considering the other process $Y_t = \sum_{i=3}^{5} X_t^{(i)}$, with the same approach as the proof of the first process, if we take the first four consecutive WV scales (i.e. $\boldsymbol{\nu}(\boldsymbol{\theta_1}) = [\nu_1^2(\boldsymbol{\theta_1}), \ldots, \nu_4^2(\boldsymbol{\theta_1})]$), we have that the determinant of the relative Jacobian matrix is $|A_{4,4}| = -{}^{2205\omega\varrho}\!/\!{}_{256}$ which implies that it is of full rank. □